\documentclass[11pt,reqno]{article}

\usepackage[margin=1in]{geometry}
\usepackage{amsmath,amsthm,amssymb}
\usepackage{graphicx}
\usepackage{bmpsize}
\usepackage{caption} 
\usepackage{esint}
\usepackage[utf8]{inputenc}

\usepackage[colorlinks=true, pdfstartview=FitV, linkcolor=blue,citecolor=blue, urlcolor=blue]{hyperref}

\usepackage[abbrev,lite,nobysame]{amsrefs}
\usepackage{times}
\usepackage[usenames,dvipsnames]{color}

\usepackage{mathtools,enumitem}

\usepackage[compact]{titlesec}

\usepackage[title]{appendix}

\usepackage{comment}

\mathtoolsset{showonlyrefs=true}
    
\newcommand{\eps}{\epsilon}

\newcommand{\grad}{\nabla}

\newcommand{\norm}[1]{\left|\left| #1 \right|\right|}
\newcommand{\abs}[1]{\left| #1 \right|}

\newcommand{\brak}[1]{\left\langle #1 \right\rangle} 

\newcommand{\dee}{\mathrm{d}}

\usepackage{bbm}

\newtheorem{theorem}{Theorem}[section]

\newtheorem*{lemma*}{Lemma}

\theoremstyle{definition}
\newtheorem{definition}[theorem]{Definition}
\newtheorem{remark}[theorem]{Remark}

\numberwithin{equation}{section}
    
\begin{document}

\title{Formalization of Landau damping in the Vlasov--Poisson equations in Lean}
\author{Jacob Bedrossian\thanks{Department of Mathematics, University of California, Los Angeles, 90095, CA USA \href{mailto:jacob@math.ucla.edu}{jacob@math.ucla.edu}. The author was supported by National Science Foundation Award DMS-2510949. The author would also like to thank Scott Armstrong for sharing some of his skill files. This work was supported in part by the OpenAI Researcher Access Program. }}

\maketitle

\begin{abstract}
We present a formalization in Lean 4 of Mouhot and Villani's theorem on nonlinear Landau damping in $\mathbb T^d$ in all Gevrey regularity indices $s > 1/3$ for small backgrounds. 
\end{abstract}

\setcounter{tocdepth}{1}
{\small\tableofcontents}

\section{Introduction}

The Vlasov equations are a fundamental model introduced by Vlasov in 1938 (see \cite{Vlasov-damping} for an English translation) for the second order motion of interacting particles in the mean-field limit. The Vlasov-Poisson equations arise when modeling charged particles interacting through electrostatic forces in plasmas and when studying matter interacting through Newtonian gravitation in galactic dynamics. We will give the precise nonlinear equations below. 

In 1946, Landau linearized the equations around the Maxwellian distribution and made the prediction of \emph{Landau damping} \cite{Landau46}, which is specifically the rapid decay (and smoothing) of the self-generated electric field, despite the equations being conservative and time-reversible. 
The rapid damping was subsequently observed in laboratory plasmas by Malmberg and Wharton \cite{MalmbergWharton64} in 1964. 
Despite the observations, physicists still raised a variety of objections to the idea that Landau damping should be mathematically present in the \emph{nonlinear collisionless} equations. Some of these objections were misguided and others were very solid\footnote{Indeed, several of the negative predictions are essentially \emph{correct} -- the regularity of the initial data determines which of the physicists' seemingly conflicting predictions is most accurate!}; see discussions in \cites{Stix,Ryutov99,MouhotVillani11} and the references therein. 
Unquestionably however, Landau damping is one of the most fundamental phenomena in the kinetic theory of plasmas (see for example \cite{Ryutov99,Stix,BoydSanderson}), and even makes an appearance as a potentially important effect in galactic dynamics \cite{Binney-Tremaine}. 

It is for these reasons, and others described below, that I believed Mouhot and Villani's theorem of Landau damping in the nonlinear Vlasov equations to be a good theorem to formalize in order to familiarize myself with agentic AI workflows for formalization using the Lean 4 proof assistant \cite{moura2021lean}. The code is available at
\begin{center}
\verb|https://github.com/Jacob24876/LandauDamping-Public.git|
\end{center}
Not only was all of the code written by agentic LLMs, but the AI was only supplied my lecture notes \cite{bedrossian2022brief} which are only a sketch that also contains several typographical errors and small mathematical inaccuracies. Other than these notes, I provided directions on how to write the two significant divergences from the source text (see below) and how to fix proofs it was stuck on, but I did not provide any tex files directly. 
Most of the work was done by Anthropic's Claude agents Sonnet 5 and Opus 4.8, with an occasional assist from OpenAI's Codex agents Sol 5.6 and Astra 6. 

Analysis is generally thought to be particularly challenging to formalize in Lean 4, especially in PDEs. 
To my knowledge, the first work formalizing a significant result in partial differential equations was the work of Armstrong and Kempe \cite{armstrong2026formalization}, which formalized the foundations of De Giorgi-Nash-Moser theory of elliptic equations. This was followed by Armstrong and Kuusi's incredible formalization of the recent quantitative stochastic homogenization work \cite{armstrong2025renormalization} (the formalization is available in \cite{armstrongKuusiCoarseGrainingRepo}). 
There had also been a variety of preceeding Lean formalizations in analysis, for example Gagliardo-Nirenberg-Sobolev inequalities \cite{van2024integrals}, Schwartz functions and tempered distributions \cite{doll2025formalizing}, and the $h$-principle \cite{van2023formalising}. There is also an ongoing larger scale project on formalizing Carleson's theorem \cite{vanDoornCarlesonRepo}. 
Weiren Sun also has a variety of introductory PDEs formalized in \cite{sunPDERepo}, and Terence Tao has formalized a significant amount of undergraduate analysis in \cite{taoAnalysisRepo}.  

Very recently of course, there was the lean formalization work of the OpenAI mathematics group for the proof of the Navier-Stokes millennium problem \cite{openaiNavierStokesEulerRepo}. However, this latter group applied an unknown --presumed incredible by academic standards-- amount of computational power to the lean formalization as well as presumed far  superior expertise at agentic coordination at large scale.

Aside from this recent work of OpenAI, to my knowledge, the present work constitutes the first formalization of a theorem regarding a nonlinear evolution equation in PDEs. In addition to the fundamental physical importance of nonlinear Landau damping and the fame of the original paper \cite{MouhotVillani11}, this result seemed a natural target for formalization since the more recent proofs are relatively non-technical and generally consist mainly of quantitative estimates on very smooth solutions instead of qualitative (i.e. there is not a lot worrying about whether you are allowed to integrate by parts and similar such things). This is largely true, although the rigorous justification of calculations in time-dependent Gevrey regularity spaces can be somewhat more delicate than in Sobolev spaces.

\subsection{Background on the Vlasov equations and Landau damping}

The Vlasov equations is a simplified model in plasma physics which neglects inter-particle collisions. 
The unknown is the \emph{distribution function} $F(t,x,v)$, which gives the number density of electrons at location $x$ moving with velocity $v$.   
If these electrons are moving due only to the electric fields generated within the plasma itself, then we obtain the model 
\begin{align*}
  & \partial_t F + v \cdot \grad_x F + E(t,x) \cdot \grad_v  F = 0 \\
  & E = \grad_x W \ast \left( \int F dv - n_0\right), 
\end{align*}
where here $n_0$ is the background ions, assumed fixed in the single species approximation\footnote{It is straightforward to extend the proof of nonlinear Landau damping to the two-species or even $n$-species model including ions of various masses and charges, however we did not bother to formalize this.}. 
The two most natural settings for this problem are $(x,v) \in \mathbb R^d \times \mathbb R^d$ and $(x,v) \in \mathbb T^d \times \mathbb R^d$, however we only consider $\mathbb T^d$ here. 
In the case of electrons in a plasma, the non-local interaction would be through Coulomb electrostatic interactions, i.e.
\begin{align}
\widehat{W}(k) = \frac{q^2}{m\eps_0 \abs{k}^2}, \label{eq:Coulomb}
\end{align}
where $q$ is the fundamental charge, $\eps_0$ is the permittivity of free space, and $m$ is the mass of an electron (the physical constants will henceforth be dropped). 
See the texts \cite{BoydSanderson,goldston2020introduction,Stix} regarding the plasma physics itself.
In the case of galactic dynamics, one only flips the sign (and changes the values of the constants) \cite{Binney-Tremaine}.

Every spatially homogeneous distribution $f^0(v)$ with $\int f^0 dv = n_0$ is a solution (notice that the electrons are all still moving around, its just that there are no density fluctuations). 
Literally the first question Vlasov asked about the equations (\cite{Vlasov-damping}) is `what happens to solutions for small disturbances of such a homogeneous equilibrium?'. 
The linearization is given by
\begin{align*}
  & \partial_t g + v \cdot \grad_x g + E \cdot \grad_v f^0 = 0 \\
  &  E = \grad_x \Delta_x^{-1} \rho \\
  & \rho(t,x) = \int g(t,x,v) dv. 
\end{align*}
The quantity $\rho$ is called the \emph{density (fluctuation)}, and it directly corresponds to ``density of electrons - average density of electrons'' as a function of $(t,x)$. 
The most important set of equilibria are 
\begin{align*}
&\textup{Maxwellians: } \quad f^0(v) = \frac{n_0}{(4\pi T)^{d/2}} e^{- \frac{\abs{v}^2}{2T}}, 
\end{align*}
for some parameters $n_0$ and $T$ (the `number density' and `temperature' respectively). 
Landau predicted rapid decay of the electric field to this linearized equation, \emph{provided the initial data had enough regularity}. 
The fact that the Vlasov equations are time-reversible means this damping must be intrinsically infinite dimensional -- the effect that shares the most in common mathematically with Landau damping is \emph{dispersion} in for example, Schr\"odinger or wave equations, \emph{not} dissipative or collisional effects like in the heat equation or Boltzmann/Landau equations. 
At least on $\mathbb T^d$, this effect can be identified as \emph{phase mixing}, wherein particles moving at different speeds tend to equi-distribute around the torus, hence smoothing out density fluctuations. 
An analogous effect in the incompressible Euler equations was discovered by Orr decades earlier in 1907 \cite{Orr07} (now called inviscid damping), however the effect of Landau damping is significantly more well known because it is a fundamental and pervasive effect in the kinetic theory of plasmas, whereas inviscid damping is harder to clearly identify as important in fluid mechanics (it is only an obvious effect in 2D while in 3D it is a rather subtle effect).

In 2008 Mouhot and Villani produced their landmark achievement \cite{MouhotVillani11}, wherein they proved that all sufficiently small analytic or Gevrey with $s$ very close to $1$ initial data leads to Landau damping solutions with dynamics essentially the same as the linearized equations.  
This proof made heavy use of Lagrangian trajectories and an explicit Newton iteration on the distribution function.
It was revolutionary and contained many brilliant insights, but was also much more complicated than is ultimately required. 
A significantly simpler proof was given by Masmoudi, Mouhot and myself in \cite{BMM13}, where we also covered down to the Gevrey threshold $s > 1/3$ conjectured in \cite{MouhotVillani11}.
The proof of Grenier, Rodnianski, and Nguyen given in \cite{GNR20,GN21} is the same general method as \cite{BMM13} but with additional improvements and simplifications. 
In my lecture notes \cite{bedrossian2022brief}, I gave a sketch of what I consider to be the simplest known proof which is essentially the proof of  Grenier, Rodnianski, and Nguyen re-written in a simpler exposition more following \cite{BMM13}, which leans on the fact that most steps in the proof can be done using classical methods from the 1990s, specifically \cite{LevermoreOliver97}. 

\section{Formalized statement}
In this section we give the exact statement formalized into lean.

Given an integer $d\geq1$ and an average-zero (real-valued) kernel $W$ on $\mathbb T^d$, we use the following
notion of solution.

\begin{definition}
Let $t_1 < t_2$. 
We say a function $f(t,x,v) \in C^1_{t,x,v}([t_1,t_2) \times \mathbb T^d \times \mathbb R^d )$ for which the following is constant in time 
\begin{align*}
n_0 := \int_{\mathbb T^d \times \mathbb R^d} f(t,x,v) \dee x \dee v 	
\end{align*}
is a classical solution to the Vlasov equations if 
\begin{align*}
\rho(t,x) = \int_{\mathbb R^d} f(t,x,v) \dee v \in C_{t,x}([t_1,t_2) \times \mathbb T^d)  	
\end{align*}
and, writing
\begin{align*}
E_f(t,x):=\grad W\ast(\rho(t,\cdot)-n_0)(x),
\end{align*}
there holds everywhere in $(t_1,t_2) \times \mathbb T^d \times \mathbb R^d$
\begin{align*}
\partial_t f + v \cdot \grad_x f + E_f \cdot \grad_v f = 0.   
\end{align*}
\end{definition}

\begin{theorem} \label{thm:MV}
Let $d\geq1$ be an integer, let $C_w>0$, and let $W$ be an average-zero, real-valued kernel on
$\mathbb T^d$ satisfying
\begin{align*}
\abs{\widehat W(k)}\leq C_w\abs{k}^{-2},\qquad k\in\mathbb Z^d\setminus\{0\}.
\end{align*}
Let $s\in(1/3,1)$, let $\lambda_0>\lambda_1>\lambda_2>0$, and let
$M_0>d/2+1$ and $m\geq 3\lfloor d/2\rfloor+7$ be integers with $m\leq M_0$.
Fix also integers $k,\ell$ with $k>d/2+1$ and $\ell>d/2$, and require
$\ell+\lfloor d/2\rfloor+2\leq m$,
specifying below the regularity class in which the solution is asserted unique. There exists
$\delta_0>0$, depending only on $C_w,d,s,M_0,m,\lambda_0,\lambda_1,$ and $\lambda_2$, with the
following property.

For every Schwartz-class function $f^0$ of $v\in\mathbb R^d$ satisfying
\begin{align*}
\sum_{\abs{\alpha}\leq M_0}
\norm{e^{\lambda_0\abs{\grad_v}}(v^\alpha f^0)}_{L^2_v}\leq\delta_0,
\end{align*}
there exist $\eps_\ast>0$ and $C>0$, depending only on the preceding fixed parameters and
on the displayed analytic bounds for $f^0$, such that the following holds. For every smooth
$h_{in}$ on $\mathbb T^d\times\mathbb R^d$, define
\begin{align*}
\eps:=\sum_{\abs{\alpha}\leq m}
\norm{e^{\lambda_1\brak{\grad_{x,v}}^s}(v^\alpha h_{in})}_{L^2_{x,v}}.
\end{align*}
If $\eps\leq\eps_\ast$ and
\begin{align*}
\int_{\mathbb T^d \times \mathbb R^d} h_{in}(x,v) \dee x \dee v = 0
\end{align*}
then there exists a function $F$ on $[0,\infty)\times\mathbb T^d\times\mathbb R^d$
whose restriction to $[0,T)$ is a classical
solution of the Vlasov equation associated to $W$, in the sense of the preceding definition, for
every $T>0$, and whose initial datum is $F(0,x,v)=f^0(v)+h_{in}(x,v)$. Moreover $F$ is the unique
such classical solution with this initial datum among those taking values continuously in
$H^k_\ell(\mathbb T^d\times\mathbb R^d)$, i.e.\ those $f(t,\cdot,\cdot)$ for which
$\sum_{\abs{\alpha}\leq \ell}\norm{v^\alpha f(t)}_{H^k_{x,v}}$ is finite and depends continuously
on $t\in[t_1,t_2)$, for every $t_1<t_2$ in the domain.

Define the perturbation $h(t,x,v):=F(t,x,v)-f^0(v)$ and its free-transport profile
\begin{align*}
g(t,z,v):=h(t,z+tv,v).
\end{align*}
There exists a smooth function $g_\infty$ of $z$ and $v$ such that, for every $t\geq0$,
\begin{itemize} 
\item  Landau damping of the electric field: 
\begin{align*}
\norm{e^{\lambda_2\brak{\grad_x,t\grad_x}^s}E_F(t)}_{L^2_x}\leq C\eps.
\end{align*}
\item Scattering to free transport:
\begin{align*}
\sum_{\abs{\alpha}\leq m-1}
\norm{e^{\lambda_2\brak{\grad_{z,v}}^s}
  (v^\alpha(g(t)-g_\infty))}_{L^2_{z,v}}
\leq C\eps e^{-\frac{\lambda_1-\lambda_2}{10}\brak{t}^s}.
\end{align*} 
\end{itemize}
\end{theorem}

\begin{remark}
There are three obvious shortcomings (and two smaller, but less obvious, short-comings). Firstly, I only treat small backgrounds. This is because the treatment of large backgrounds would require formalizing a great deal of complex analysis theory surrounding Laplace transforms and this was deemed not worth the investment. 
The other four shortcomings were much smaller cut-corners to complete the formalization faster: the requirement $ s < 1$ (a restriction on the decay rate, not the initial data), the unnecessarily large requirements on the velocity moments, the Schwartz class requirement on $f^0$ (a requirement on having faster than polynomial decay in velocity), and the fact that the theorem is technically only stated for forward time, while by time-reversibility, it must also hold backward in time.  
\end{remark}

The proof itself follows the lecture notes \cite{bedrossian2022brief} and the original sources \cite{BMM13,GNR20} closely except for two departures. 
The simplest departure is the treatment of the Volterra equation. In \cite{BMM13} this is treated highly inefficiently, a weakness that is corrected by the much more elegant and natural approach taken in \cite{GNR20}. However, these approaches use the entire machinery of Laplace transforms, which would have been a tedious job to formalize. Hence, our formalization only includes a straightforward fixed point argument that covers all small backgrounds. The significant departure is the need to develop a well-posedness theory and to justify the qualitative calculations necessary to do the energy estimates.  
This is a slightly delicate task to do in Gevrey regularity due to the necessity of having constantly decreasing regularity index. 
The approach we took is to regularize the nonlinearity and use a Picard theorem for existence (essentially as in classical quasilinear existence theories in PDEs such as in the proof of local well-posedness of the Euler equations in \cite{MajdaBertozzi}), derive all of the quantitative estimates there with regularity to spare, and then pass to the limit uniformly (although it is more delicate here due to Gevrey regularity and having to control velocity moments). 
We did not bother to formalize a full local well-posedness theory for solutions of arbitrary distance from equilibrium (except for the uniqueness assertion), however, we do not believe this would have been particularly difficult.    

\section{AI Workflow} 

This was a learning project to force me to learn how to use agentic AI. 
As such, the workflow evolved several times as I learned lessons the hard way.

As remarked above, the AIs were only directly supplied my lecture notes \cite{bedrossian2022brief}. The details depended on the workflow, but the AI would first draft the requisite detailed tex files (around 5-10 pages each), figure out what `infrastructure' was required (for example, Fourier analysis on cylinder domains) and then draft those tex files. Then, it would make a formalization plan, then do the actual formalization, then audit.
The work basically evolved in three rough stages. The first workflow was experimental and disorganized and not suitable for large scale lean. The second workflow was an over-engineered semi-automatic formalization that involved the orchestration of 5 separate agents. When this worked it worked fully automatic. However, it sometimes felt that more time was spent on writing and debugging bureaucratic book-keeping than on actual formalization, and it struggled to react effectively when it hit a snag. 
The last workflow was the only useful one: (a) have the agents fully work out and verify (with symbolic checking tools such as SymPy) the proofs in natural language first as one would expect, and then pass to a streamlined set of 3 agents (planner, formalizer, and auditor) to convert into lean and audit for faithfulness.  

Here is a short summary of advice to a mathematician trying to formalize on a budget:  
\begin{itemize}
\item[(i)] You should have a more solid and complete source than what I used, which will help the agents keep on track far better (you can still have the AI write it, but it should be done first, preferably with an explicit directed graph of mathematical dependence pre-made). Regardless of exactly how you organize the source, note that agents struggle to reason directly in lean. 
Even if they hit an unexpected snag mid-formalization, make your agents reason and write proofs in natural language and use a real prover-checker loop with symbolic tools such as SymPy.  
\item[(ii)] Have your agents first build a lean scaffolding from your source, i.e. lean statements with empty proofs of the final statement and the main lemmas/propositions.  
\item[(iii)] Have your agents pull Scott Armstrong's repositories and glean good ideas, including how to formalize $L^p$ and Sobolev spaces efficiently, how faithfulness is ensured, and how final lean statements are audited.   
\item[(iv)] Have your agents plan formalization tasks into clear blueprints before executing, just like you would when working on a large-scale software engineering task. Do not try to formalize too much at once.    
\item[(v)] Have your agents periodically reflect on how to improve their workflow in terms of time, elaboration costs, correctness/faithfulness, and AI usage costs. Have them re-write their behavior accordingly. This is the easy version of various strategies used to enhance the effectiveness of LLM agents by updating plain text memories based on self-reflection; see e.g. \cite{shinn2023reflexion,zhao2024expel,wang2024agent}.  
\end{itemize}
I experimented with various ways to cut usage such as maintaining an index file system and various approaches to orchestration but I am not confident these were productive, and some of my orchestration attempts were actively unproductive.

\bibliographystyle{abbrv}
\bibliography{eulereqns,JacobBib,VladBib}

\def\cprime{$'$} \def\cprime{$'$}
\begin{bibdiv}
\begin{biblist}

\bib{armstrong2026formalization}{article}{
      author={Armstrong, Scott},
      author={Kempe, Julia},
       title={Formalization of de giorgi--nash--moser theory in lean},
        date={2026},
     journal={arXiv preprint arXiv:2604.05984},
}

\bib{armstrong2025renormalization}{article}{
      author={Armstrong, Scott},
      author={Kuusi, Tuomo},
       title={Renormalization group and elliptic homogenization in high contrast},
        date={2025},
     journal={Inventiones mathematicae},
      volume={242},
      number={3},
       pages={895\ndash 1086},
}

\bib{armstrongKuusiCoarseGrainingRepo}{misc}{
      author={Armstrong, Scott},
      author={Kuusi, Tuomo},
       title={Coarsegraining: A lean 4 formalization of coarse-graining theory for elliptic equations},
        date={2026},
        note={\url{https://github.com/scottnarmstrong/CoarseGraining}},
}

\bib{bedrossian2022brief}{article}{
      author={Bedrossian, Jacob},
       title={A brief introduction to the mathematics of landau damping},
        date={2022},
     journal={arXiv preprint arXiv:2211.13707},
}

\bib{BMM13}{article}{
      author={Bedrossian, Jacob},
      author={Masmoudi, Nader},
      author={Mouhot, Clement},
       title={Landau damping: paraproducts and gevrey regularity},
        date={2016},
     journal={Annals of PDE},
      volume={2},
      number={1},
       pages={1\ndash 71},
}

\bib{Binney-Tremaine}{book}{
      author={Binney, J.},
      author={Tremaine, S.},
       title={Galactic dynamics},
   publisher={Princeton University Press (2d edition)},
        date={2008},
}

\bib{BoydSanderson}{book}{
      author={Boyd, T. J.~M.},
      author={Sanderson, J.~J.},
       title={The physics of plasmas},
   publisher={Cambridge University Press},
     address={Cambridge},
        date={2003},
        ISBN={0-521-45290-2; 0-521-45912-5},
         url={http://dx.doi.org/10.1017/CBO9780511755750},
      review={\MR{1960956 (2005e:82113)}},
}

\bib{doll2025formalizing}{article}{
      author={Doll, Moritz},
       title={Formalizing schwartz functions and tempered distributions},
        date={2025},
     journal={arXiv preprint arXiv:2510.24060},
}

\bib{goldston2020introduction}{book}{
      author={Goldston, Robert~J},
       title={Introduction to plasma physics},
   publisher={CRC Press},
        date={2020},
}

\bib{GN21}{article}{
      author={Grenier, Emmanuel},
      author={Nguyen, Toan},
       title={Generator functions and their applications},
        date={2021},
     journal={Proceedings of the American Mathematical Society, Series B},
      volume={8},
      number={20},
       pages={245\ndash 251},
}

\bib{GNR20}{article}{
      author={Grenier, Emmanuel},
      author={Nguyen, Toan~T},
      author={Rodnianski, Igor},
       title={Landau damping for analytic and gevrey data},
        date={2020},
     journal={arXiv preprint arXiv:2004.05979},
}

\bib{Landau46}{article}{
      author={Landau, Lev},
       title={On the vibration of the electronic plasma},
        date={1946},
     journal={J. Phys. USSR},
      volume={10},
      number={25},
}

\bib{LevermoreOliver97}{article}{
      author={Levermore, D.},
      author={Oliver, M.},
       title={Analyticity of solutions for a generalized {Euler} equation},
        date={1997},
     journal={J. Diff. Eqns.},
      volume={133},
       pages={321\ndash 339},
}

\bib{MajdaBertozzi}{book}{
      author={Majda, A.},
      author={Bertozzi, A.~L.},
       title={Vorticity and incompressible flow},
   publisher={Cambridge University Press},
        date={2002},
}

\bib{MalmbergWharton64}{article}{
      author={Malmberg, J.},
      author={Wharton, C.},
       title={Collisionless damping of electrostatic plasma waves},
        date={1964},
     journal={Phys. Rev. Lett.},
      volume={13},
      number={6},
       pages={184\ndash 186},
}

\bib{MouhotVillani11}{article}{
      author={Mouhot, Cl{\'{e}}ment},
      author={Villani, C{\'{e}}dric},
       title={On {Landau} damping},
        date={2011},
     journal={Acta Math.},
      volume={207},
       pages={29\ndash 201},
}

\bib{moura2021lean}{inproceedings}{
      author={Moura, Leonardo~de},
      author={Ullrich, Sebastian},
       title={The lean 4 theorem prover and programming language},
organization={Springer},
        date={2021},
   booktitle={International conference on automated deduction},
       pages={625\ndash 635},
}

\bib{openaiNavierStokesEulerRepo}{misc}{
      author={{OpenAI}},
       title={Lean certificates accompanying {Navier--Stokes} and {Euler} results},
        date={2026},
        note={\url{https://github.com/openai/NavierStokesAndEuler}},
}

\bib{Orr07}{article}{
      author={Orr, W.},
       title={The stability or instability of steady motions of a perfect liquid and of a viscous liquid, {Part I}: a perfect liquid},
        date={1907},
     journal={Proc. Royal Irish Acad. Sec. A: Math. Phys. Sci.},
      volume={27},
       pages={9\ndash 68},
}

\bib{Ryutov99}{article}{
      author={Ryutov, DD},
       title={Landau damping: half a century with the great discovery},
        date={1999},
     journal={Plasma physics and controlled fusion},
      volume={41},
      number={3A},
       pages={A1},
}

\bib{shinn2023reflexion}{article}{
      author={Shinn, Noah},
      author={Cassano, Federico},
      author={Gopinath, Ashwin},
      author={Narasimhan, Karthik},
      author={Yao, Shunyu},
       title={Reflexion: Language agents with verbal reinforcement learning},
        date={2023},
     journal={Advances in neural information processing systems},
      volume={36},
       pages={8634\ndash 8652},
}

\bib{Stix}{book}{
      author={Stix, T.},
       title={Waves in plasmas},
   publisher={Springer},
        date={1992},
}

\bib{sunPDERepo}{misc}{
      author={Sun, Weiran},
       title={{PDE Lean Formalization}},
        date={2026},
        note={\url{https://github.com/weiran-sun/pde}},
}

\bib{taoAnalysisRepo}{misc}{
      author={Tao, Terence},
       title={A lean companion to analysis i},
        date={2026},
        note={\url{https://github.com/teorth/analysis}},
}

\bib{vanDoornCarlesonRepo}{misc}{
      author={van Doorn, Floris},
       title={Formalization of a generalized carleson's theorem},
        date={2026},
        note={\url{https://github.com/fpvandoorn/carleson}},
}

\bib{van2024integrals}{inproceedings}{
      author={van Doorn, Floris},
      author={Macbeth, Heather},
       title={Integrals within integrals: A formalization of the gagliardo-nirenberg-sobolev inequality},
organization={Schloss Dagstuhl--Leibniz-Zentrum f{\"u}r Informatik},
        date={2024},
   booktitle={15th international conference on interactive theorem proving (itp 2024)},
       pages={37\ndash 1},
}

\bib{van2023formalising}{inproceedings}{
      author={van Doorn, Floris},
      author={Massot, Patrick},
      author={Nash, Oliver},
       title={Formalising the h-principle and sphere eversion},
        date={2023},
   booktitle={Proceedings of the 12th acm sigplan international conference on certified programs and proofs},
       pages={121\ndash 134},
}

\bib{Vlasov-damping}{article}{
      author={Vlasov, A.~A.},
       title={The vibrational properties of an electron gas},
        date={1938},
     journal={Zh. Eksp. Teor. Fiz.},
      volume={291},
      number={8},
        note={In russian, translation in english in {\em Soviet Physics Uspekhi}, vol. 93 Nos. 3 and 4, 1968},
}

\bib{wang2024agent}{article}{
      author={Wang, Zora~Zhiruo},
      author={Mao, Jiayuan},
      author={Fried, Daniel},
      author={Neubig, Graham},
       title={Agent workflow memory},
        date={2024},
     journal={arXiv preprint arXiv:2409.07429},
}

\bib{zhao2024expel}{inproceedings}{
      author={Zhao, Andrew},
      author={Huang, Daniel},
      author={Xu, Quentin},
      author={Lin, Matthieu},
      author={Liu, Yong-Jin},
      author={Huang, Gao},
       title={Expel: Llm agents are experiential learners},
        date={2024},
   booktitle={Proceedings of the aaai conference on artificial intelligence},
      volume={38},
       pages={19632\ndash 19642},
}

\end{biblist}
\end{bibdiv}
\end{document}